\documentclass[reqno,11pt,a4paper,oneside]{article}
\usepackage{setspace}
\usepackage{amsfonts}
\usepackage{amsthm}
\usepackage[mathscr]{euscript}
\usepackage[cp1250]{inputenc}
\usepackage[T1]{fontenc}
\usepackage{amssymb}
\usepackage{amsmath}
\usepackage[dvipdf]{graphicx}
\def\={\discretionary{-}{-}{-}}

\newtheorem{tw}{Theorem}

\newtheorem{lm}{Lemma}
\newtheorem{rem}{Remark}

\newtheorem{df}{Definition}
\newtheorem{wn}{Corollary}
\newtheorem{prz}{Example}

\title{Weierstrass functions and a generalization of the  additive-multiplicative Weierstrass  inequality }
\author{Halina Wi\'sniewska }
\date{Institute of Mathematics,  Kazimierz Wielki University, 85-072 Bydgoszcz, Poland\\ email: halinkaw@ukw.edu.pl}
\begin{document}
\doublespacing
\maketitle

\vspace{12pt}

\begin{abstract}
 Let $J$ denote the interval either $(0,1]$ or $ [1, \infty)$. A positive function $f$ on $J$ with $f(1) =1$ is reffered to as a Weierstrass function if it fulfils the double inequality for $x,y \in J$:
 \begin{equation} \label{q}
 f(x) + f(y) -1 \leq f(xy) \leq f(x)f(y).
 \end{equation}
 By means of such functions we can extend the classical Weierstrass inequality (inequality (\ref{q}) for $f(x) = x$) to some trigonometric, Euler gamma, and log- functions.
 Utilizing the Weierstrass property of $f(x) = \frac{\ln (1+x)}{\ln2}$,  we obtain a new multiplicative inequality which, in turn, generalizes the classical Weierstrass inequality.

\end{abstract}


\vspace{12pt}


\section{Introduction} \label{wprowadzenie}
Let $a_i \in (0,1)$ for $i=1,2,...,n$, $n \geq 2$. The classical Weierstrass's inequality asserts that:
\begin{equation} \label{klasycz}
\prod_{i=1}^{n} (1-a_i) \geq 1 - \sum_{i=1}^{n} a_i.
\end{equation}
\cite[p. 210]{Mitrinovic}.
Inequality (\ref{klasycz}) is one of the most important inequalities in analysis, and  it was the object of interest of many researchers who extended (\ref{klasycz}) by means of known techniques  (see \cite{Klamkin4, Klamkin5,Nizar,Pec,WU} and the references therein).

\vspace{12pt}

The goal of this paper is an extension of the Weierstrass inequality in a different way, as follows.

\vspace{12pt}

\noindent
By substituting in inequality (\ref{klasycz}): $x_i = 1-a_i $, $i=1,2,, ..., n \geq 2$,   we obtain an equivalent form of (\ref{klasycz}):
\begin{equation} \label{klasycz-po-podst}
\prod_{i=1}^{n} x_i \geq \sum_{i=1}^{n} x_i - (n-1), \ with \ x_i \in (0,1], \ n \geq 2.
\end{equation}
To simplify  our calculations,  we  consider an inequality equivalent to  (\ref{klasycz-po-podst}) for $n=2$ only:
\begin{equation} \label{W}
x+y -1 \leq xy, \ with \ x,y \in (0,1].
\end{equation}
Let us note that 
\begin{equation} \label{nr}
 inequality \  (\ref{W}) \  also \ holds\  for \  x,y \in [1, \infty) \  (\mbox{see} \  \cite[\rm{Theorem \ 5}]{Finol-Wojt}).
\end{equation}

 The main goal of this paper is to extend the above inequalities (\ref{W}) and (\ref{nr}) by replacing $x$ and $y$ by the values $f(x)$ and $f(y)$ of a certain real function $f$,  and then  inserting the value $f(xy)$  between such modified inequalities:

\begin{equation} \label{def-nier-Weier}
f(x) + f(y) -1 \leq f(xy) \leq f(x)f(y) \ with \  x,y \in (0,1]\  or \ x,y \in [1, \infty).
\end{equation}

 A precise definition is given below.

\begin{df} \label{d-fun-Weiers}
Let $J$ be one of the two intervals: either $(0,1]$ or $[1, \infty)$, and let  $f$ be a positive function on  $J$ with $f(1)=1$. The function $f$ is said to be a Weierstrass function on $J$, or that $f$ has the Weierstrass property on $J$, if it fulfils condition $(\ref{def-nier-Weier})$;  equivalently, by induction on $n$:

\begin{equation} \label{prodWeier}
\sum_{i=1}^{n} f(x_i) - (n-1) \leq f \left ( \prod_{i=1}^{n} x_i \right ) \leq \prod_{i=1}^{n} f(x_i), \ with \ x_1, ..., x_n \in J.
\end{equation}
\end{df}

From the definition  it follows that every Weierstrass function on  $J$ is submultiplicative. Such functions were studied in many papers (see e.g. \cite{Alzer1,Finol-Wojt, Maligranda} and the references therein), and there are known conditions for a given function on $J$ to be submultiplicative. In our paper we deal with relationships between the values $f(xy)$ and $f(x) + f(y) -1$, with $x,y \in J$. For this purpose, we define $l-$ and $r-$Weierstrass functions.

\begin{df} \label{def-l-r-Weierstrass}
Let $J$ be  either $(0,1]$ or $[1, \infty)$, and let  $f$ be a positive function on  $J$ with $f(1)=1$. The function $f$ is said to be left-Weierstrass $(l-$Weierstrass$)$ if it  fulfils the left side of inequality $(\ref{def-nier-Weier})$:
\begin{equation} \label{rn-1}
f(xy) \geq f(x) + f(y) -1,
\end{equation}
or, equivalently, the left side of inequality $(\ref{prodWeier})$.

If $f$ fulfils the reversed inequality to $(\ref{rn-1})$ on $J$ then $f$ is referred to as a right-Weierstrass $(r-$Weierstrass$)$ function.
\end{df}

\begin{rem} \label{rem11}
\rm{Notice that if} $f: J \to J$ is a $r-$Weierstrass function then it is already a submultiplicative function $($by $(\ref{W})$ or $(\ref{nr}))$.
\end{rem}

\noindent
The notion of a Weierstrass function allows us to extend  immediately inequalities (\ref{W}) and (\ref{nr}) as follows:

\begin{tw} \label{t1-odwr}
Let $J$ denote one of the two intervals: either $(0,1]$ or $[1, \infty)$, and let $f$ be an automorphism of  $J$ with $\Psi$ the inverse function of $f$. If $f$ is a  Weierstrass function, then:
\begin{equation} \label{dla-psi}
x+y -1 \leq f(\Psi (x) \Psi (y)) \leq xy, \ with \ x,y \in J;
\end{equation}
inductively, for $n \geq 2$ and $x_1, ..., x_n \in J$,
\begin{equation}
\sum_{i=1}^{n} x_i - (n-1) \leq f \left ( \prod_{i=1}^{n} \Psi (x_i) \right ) \leq \prod_{i=1}^{n} x_i.
\end{equation}
\end{tw}

Of course, the identity function  $f(x) = x$ is a Weierstrass function both   on $(0,1]$ and  $[1, \infty)$. Thus, a natural questions arise:
\begin{itemize}
\item[$(Q1)$] Are there any other Weierstrass functions  on  either  $(0,1]$ and  $[1, \infty)$ then the identity function?
\item[$(Q2)$] What conditions must satisfy a real function to be a Weierstrass function?
\item[$(Q3)$] Is every submultiplicative function a Weierstrass function?
\end{itemize}
In this paper, we  give positive answers to  questions $(Q1)$ and $(Q2)$ to functions of class $C^1$: in the next section, we  give a sufficient condition for a given real function $f$ to be a Weierstrass function; it applies  well to some classical functions, they products, and compositions. In Example \ref{Ex4} we give a negative answear to question $(Q3)$.

\section{Main Criteria} \label{Main-Results}

In the theorem below, $f$ denotes a function defined either on  $J= (0,1]$ or on $J =[1.\infty)$.

\begin{tw} \label{o-monot}
Let  $f$  be a positive function on $J$ of class $C^1$ with $f(1) =1$.  
A sufficient condition for $f$ to be an $l-$Weierstrass function is this:
\begin{equation} \label{star}
\mbox{The function} \ H_f(x) := x f'(x) \  \mbox{is non-decreasing on}  \ J. 
\end{equation}
If $H_f$ is strictly increasing, then  inequality $(\ref{rn-1})$ is strict.

These statements remain valid on replacing everywhere '$l-$Weierstrass' by '$r-$Weierstrass' and '$[$non-$]$decreasing' by '$[$non-$]$increasing'.
\end{tw}

\noindent
\textbf{Proof of Theorem \ref{o-monot}. }

We shall prove only the case $f$ is an $l-$Weierstrass function on $J= (0,1]$ because the remaining cases can be proved in a similar fashion.

Let   $y$ be a fixed arbitrary element of $J$, and set 
\begin{equation} \label{phi}
\varphi (x) := f(xy) - f(x) - f(y) +1, \ x \in J.
\end{equation}
Because $\varphi (1) =0$, to prove that $f$ is $l-$Weierstrass it is enough to show that $\varphi$ is non-increasing on $J$.

\noindent
Since $f$ is of class $C^1$, this is equivalent to the inequality
\begin{equation} \label{n2}
\varphi '(x) \geq 0, \ with \ x \in J, \ i.e.,
\end{equation}

\begin{equation} \label{no1}
 yf' (xy) \leq f'(x), \ with \ x \in (0,1).
\end{equation}
 Multiplying both sides of   (\ref{no1}) by $x$ and setting $s(y) : = xy < x$ we obtain that (\ref{no1}) is equivalent to the condition:

\begin{equation} \label{no2}
s(y) \cdot f' (s(y)) \leq x \cdot f' (x),  \ with \  x \in (0,1], 
\end{equation}
Since $y \in J$ was arbitrary fixed, from (\ref{no2}) we obtain that
\begin{equation} \label{no3}
 s \cdot f' (s) \leq x f'(x), \ whenever \ s<x, \ with \ x \in (0,1), \ i.e.,
\end{equation}
 the function $H_f (x) = x \cdot f' (x)$ is non-decreasing on $(0,1)$.

On the other hand, every $y \in J$ is of the form $y = \frac{s}{x}$ for some $s,x \in (0,1)$; thus, from (\ref{no3}) we obtain (\ref{no1}). 

We thus have proved that $f$ is an $l-$Weierstrass function on $J$ provided that $H_f$ is non-decreasing  on $J$. It is also obvious that the inequality (\ref{rn-1}) is strict if $H_f$ is strictly increasing on $J$. $\square$

\vspace{12pt}

For functions of class $C^2$ we obtain a more usefull criterion then in Theorem \ref{o-monot}.
The proof is obvious.

\begin{wn} \label{C-2}
Let $f$ be a positive function of class $C^2$ with $f(1) =1$. Then $f$ has the $l-$ $[r-]$Weierstrass property if it  fulfils the condition:
\begin{equation} \label{G}
\mbox{The  expression} \ G_f (x) : = f' (x) \cdot \left ( \frac{1}{x} + \frac{f'' (x)}{f' (x) }\right ) \ \mbox{is non-negative} \left [ \mbox{non-positive} \right ] on \  J.
\end{equation}
\end{wn}

From Theorem \ref{o-monot} we immediatelly obtain:

\begin{wn}  \label{lem1}
Let $f_1, ..., f_k$, $k \geq 2$ be submultiplicative functions on $J$ of class $C^1$ with $f_j (1) =1$, $j =1, ...,k$ satisfying condition (\ref{star}). Then  $F_k := f_1 \cdot f_2 \cdot ... \cdot f_k$ is a  Weierstrass function on $J$.
\end{wn}

\noindent
\textbf{Proof of Corollary \ref{lem1} }

It follows from the formula:
$
x {(f (x) \cdot g (x) )}' = x f' (x) g (x) + x f (x)g' (x),
$
 condition (\ref{star}) and the assumptions of the corollary. $\square$

\begin{rem}
 \rm{If a function} $f$ of class $C^1$ is submultiplicative on $J$ with $f(1) =1$, and satisfies condition $(\ref{star})$ then, for every $ k \geq 2$, the function $f^k$  is  Weierstrass, too, but in Corollary \ref{dz} below this is a particular case of a more general result.
\end{rem}

\begin{wn} \label{lem2}
Let $f$ and $g$ be two functions  $J \to J$ of class $C^1$, with $f(1) = g(1) =1$, and with $f,g$ and $g'$ non-decreasing. If $f$ satisfies condition $(\ref{star})$, then the superposition $g \circ f$ satisfies condition $(\ref{star})$, too, and hence the function $g \circ f$ is $l-$Weierstrass.

 Additionally, if  $f$ and $g$ are submultiplicative on $J$, then the superposition $g \circ f$ is a Weierstrass function.
 \end{wn}
 
 \noindent
\textbf{Proof of Corollary \ref{lem2} }

It follows from the formula:
$
x {g(f(x))}' = x g' (f(x)) f' (x),
$
 condition (\ref{star}) and the assumptions of the corollary. $\square$
 
Because the function $g(x) = x^{\alpha}$ is increasing, multiplicative on $[0, \infty)$ with $g$ convex iff $\alpha \geq 1$, we immediatelly obtain:
 
\begin{wn} \label{dz} 
 If $f$ is a Weierstrass function of class $C^1$ on $J$ then, for every $\alpha > 1$, the function $F_{\alpha}  := f^{\alpha} $ is a Weierstrass function, too.
\end{wn}

\section{Examples of  Weierstrass functions}
\noindent
We present below   some examples of Weierstrass functions on proper intervals.

\begin{prz} \label{prz3}
 \rm{By} \cite[Theorem 12 $(iv)$]{Finol-Wojt}, the function $f(x) = \frac{4}{\pi} \arctan x$ is submultiplicative on $J =(0,1]$, with $f(1) = 1$. 
 Because the function $H_{f} (x) = \frac{4}{\pi} \cdot \frac{x}{1+ x^2}$ is increasing on $J$, from Theorem \ref{o-monot} we obtain that $f$ has the $l-$Weierstrass property on $J$;
 hence by Theorem \ref{t1-odwr},
 $$
x+y -1 \leq \frac{4}{\pi} \arctan \left ( \tan \frac{\pi}{4} x \cdot \tan \frac{\pi}{4} y \right )  \leq xy, \ with \ x,y \in (0,1].
$$

\end{prz}

\begin{prz} \label{prz2}
 \rm{By the proof similar to that of} \cite[Theorem 12 $(i)$]{Finol-Wojt}, the function $ f(x) = \sin \frac{\pi}{4} x$ is submultiplicative on $J=(0,1]$, with $f(1) =1$.  Because the function $H_f(x) = xf' (x) = \frac{\pi x}{4} \cos \frac{\pi}{4}x$ is increasing on $J$,  from Theorem \ref{o-monot} it follows that  $f$ is a Weierstrass  function on $J$. Since  $\Psi (x) = f^{-1} (x) = \frac{4}{\pi} \arcsin x$, from Theorem \ref{t1-odwr} we obtain that: 
$$
x+y -1 \leq \sin \left ( \frac{4}{\pi} \arcsin x \cdot \arcsin y \right ) \leq xy, \ with \ x,y \in (0,1]. 
$$
\end{prz}

As we have noticed in Remark \ref{rem11}, every Weierstrass function is already submultiplikative (by definition).
In the next example we show that the class of submultiplicative functions is broader then the class of a Weierstrass functions.

\begin{prz} \label{Ex4}
 \rm{From the proof of}  \cite[Theorem 2]{Finol-Wojt}   it follows that the function $f(x) = \frac{\cos x}{\cos 1}$ is submultiplicative on $J =(0,1]$, with $f(1) =1$. Here $H_f (x) = \frac{-x \sin x}{\cos 1}$ is a decreasing function on $J$. By Theorem \ref{o-monot},  $f$ is not a Weierstrass function.
\end{prz}

\section{The case of  a  log function}

Let, for $k >0$, the symbol $\log_k x$ denote the logarithm of $x > 0$ in base $k$.
In this section, we deal with the Weierstrass property of the function $L(x) = \log_2 (1+x)$ and its consequences.

It is known that 
\begin{equation} 
\label{kk}
L \  \mbox{\textit{is submultiplicative both on}} \  (0,1] \  \mbox{\textit{and}} \  [1, \infty),  \ \mbox{\textit{separately}};
\end{equation}
(see \cite[Theorem 10 (i)]{Finol-Wojt}).
In the theorem below, we complete this property into the direction of the $l-$Weierstrass property for $L$.

\begin{tw} \label{tw-dla-logar}
Let $J$ denote either $(0,1]$ or $[1, \infty)$.  Then:
\begin{itemize}
\item[$(i)$] The function $L$ has the $l-$Weierstrass property on $J$, and the inequality is strict for $x \neq y$ with $x,y \in J$.
Hence, $L$ is a Weierstrass function on $J$.
\item[$(ii)$] For every $n \geq 2$, and $x_1, ..., x_n \in J$,
\begin{equation} \label{gw1}
\prod_{i=1}^n (1+ x_i) \leq 2^{n-1} \cdot \left ( 1+ \prod_{i=1}^n x_i \right ).
\end{equation}
\end{itemize}
\end{tw}

\begin{rem}
\rm{Inequality} $(\ref{gw1})$ seems to be new. It is also another generalization of the Weierstrass inequalities $(\ref{W})$ and $(\ref{nr})$.
Indeed, after multiplying and rearranging of numbers in (\ref{gw1}),

\noindent
for $n =2$ we obtain $(\ref{W})$ and $(\ref{nr})$; 

\noindent
for
 $n =3$ we obtain $$x_1 (1+x_2) + x_2 (1+x_3) + x_3 (1+x_1) \leq 3 (1 + x_1 x_2 x_3);$$
and so on.
\end{rem}

\noindent
\textbf{Proof of Theorem \ref{tw-dla-logar}}

\noindent
$(i)$ For arbitrary $x \in (0, \infty)$, we have $H_f (x) = \frac{x}{(1+x) \ln 2}$, thus $H_f$ is strictly increasing both on $J = (0,1]$ and on $J = [1, \infty)$. Now we apply Theorem \ref{o-monot} and (\ref{kk}).

\noindent
$(ii)$ By part $(i)$ of our theorem and Definition \ref{def-l-r-Weierstrass}, 
$$
\log_2 \left ( \prod_{i=1}^n (1+ x_i) \right ) = \sum_{i=1}^n \log_2 (1+ x_i) \leq (n-1) + \log_2 \left ( 1 + \prod_{i=1}^n x_i \right ) =
$$
\begin{equation} \label{mm1}
 \log_2 \left ( 2^{n-1} \left ( 1 + \prod_{i=1}^n x_i \right ) \right ), \ x_i \in J,
\end{equation}
 whence inequality (\ref{gw1}) follows. $\square$

\section{The case of the Euler gamma function}

Let $\Gamma_{\{a\}} (x) = \frac{\Gamma (ax)}{\Gamma (a) }$, for $a >0$ fixed and $x > 0$ arbitrary. In this section, we study the Weierstrass property for functions of the form   $\Gamma_{\{a\}}$. It is already known that \cite[Theorem  8 (i)]{Finol-Wojt}:

\begin{lm}
The function $\Gamma_{\{a\}}$ is submultiplicative on $J = (0,1]$ for $a \in (0, \xi]$, where $\xi  = 0,21609...$ is the solution of the equation $\sum_{n=1}^{\infty} \frac{2nx + x^2}{n {(n + x)}^2} = \gamma,$ where $\gamma$ is the Euler's constant $0,577215...$.
\end{lm}
\noindent
 For our purposes, we complete the lemma essentially, as follows:
let $x_{min} = 1,4616...$ denote the point  where  $\Gamma$ reaches its only minimum on $(0, \infty)$ (see \cite[p. 303]{Pier}), and set $x_1 = x_{min} -1 = 0,4616...$; then: 

\begin{tw} \label{tw-o-gamma}
With the notation as above, the function $\Gamma_{\{a\}}$ has the $l-$Weierstrass property on $J = (0,1]$ for every $a \in (0, x_1]$.

In particular, $\Gamma_{\{a\}}$ is a Weierstrass function on $J$ for every $a \in (0, \xi]$.
\end{tw}

The proof of our theorem  is the result of log-convexity of the function $\Gamma$ on $(0, \infty)$ (cf. \cite[Theorem 7.71]{Stromberg}) and of Corollary \ref{C-2}. 

Let us recall that a function $F$ on an interval $(a,b)$ is log-convex if the function $\ln F$ is convex.  It is well known that, if $F$ is of class $C^2$ on $(a,b)$ then $F$ is log-convex iff
\begin{equation} \label{D}
F'' (x) \cdot F (x) - {(F' (x))}^2 \geq 0 \ on \ (a,b).
\end{equation}

In the proof of our theorem we apply the following lemma:

\begin{lm} \label{lem-dwoch-fun}
If $f$  in Theorem \ref{o-monot} is of class $C^2$ and log-convex, then $f$ has the $l-$Weierstrass property in two either cases:
\begin{itemize}
\item[$(i)$] The both functions, $f(x)$ and $x f(x)$, are strictly decreasing on $J$,
\item[$(ii)$]  The both functions, $f(x)$ and $xf(x)$, are strictly increasing on $J$.
\end{itemize}
\end{lm}

\noindent
\textbf{Proof of Lemma \ref{lem-dwoch-fun}.}
We  consider only case $(i)$ because the proof of case $(ii)$ is similar.

We shall show that the function 
\begin{equation} \label{krop-1}
G_f  (x) =  f'(x) \left (  \frac{1}{x} + \frac{f'' (x)}{f' (x)} \right )
\end{equation}
in Theorem \ref{tw-dla-logar} is positive on $J$.

 Since, by hypoteses,
\begin{equation} \label{krop-2}
f(x) \cdot f'' (x) - {(f' (x))}^2  \geq 0,
\end{equation}
and $f' (x) < 0$ on $J$, inequality (\ref{krop-2}) is equivalent to 
$$
\frac{f'' (x)}{f' (x)} \leq \frac{f' (x)}{f(x)},
$$
whence
$$
\frac{1}{x} + \frac{f'' (x)}{f' (x)} \leq {(\ln (x \cdot f(x))}'.
$$
Then, in (\ref{krop-1}), we have
\begin{equation} \label{krop-3}
G_f  (x) \geq    f' (x) \cdot {(\ln (x \cdot f(x))}' .
\end{equation}
Because, by hypoteses, the functions $f(x)$ and $xf(x)$ are decreasing simultaneously on $J$, the right side of (\ref{krop-2}) is positive on $J$. Now we apply Corollary  \ref{C-2}. $\square$

\noindent
\textbf{Proof of Theorem \ref{tw-o-gamma}.}
Notice first that $\Gamma_{\{a\}}$ is log-convex on $(0, \infty)$ (this follows from the log-convexity of $\Gamma$ 
 \cite[Theorem 7.71 and Exercise 8 (g), p.472]{Stromberg} and inequality (\ref{D}) applied for $F =\Gamma_{\{a\}}$ ), and hence we can apply Lemma \ref{lem-dwoch-fun}.

 For $x \in (0, 1]$ we have the identity:
 $$
 x \cdot \Gamma_{\{a\}} (x) = \frac{ax \Gamma (ax)}{a \Gamma (a)} = \frac{\Gamma (ax + 1)}{\Gamma (a+1)},
 $$
and hence  the function $x \cdot \Gamma_{\{a\}} (x)$ is decreasing on $J= (0,1]$ for $ax + 1 \leq x_{min}$, i.e., for every $a \leq x_{min} -1 = x_1$.

Additionally, the function $\Gamma_{\{a\}}$ is decreasing on $J$ for $0 < a \leq x_{min}$. Summing up, case $(i)$ of Lemma \ref{lem-dwoch-fun} holds true for $a \leq x_1$.
$\square$

From Theorem \ref{tw-o-gamma} and Definition \ref{def-l-r-Weierstrass}, we immediately obtain  new multiplicative-additive inequalities of $\Gamma$.

\begin{wn}
For every $a < x_1$, and every $x,y \in  (0,1]$, we have
$$
\Gamma_{\{a\}} (xy) \geq \Gamma_{\{a\}} (x) + \Gamma_{\{a\}} (y) -1.
$$
Equivalently:
\begin{equation} \label{e1}
\Gamma (axy) \geq \Gamma (ax) + \Gamma (ay) - \Gamma (a).
\end{equation}
\end{wn}

Setting in $(\ref{e1})$: $u = ax$, and $v = ay$ with $x,y \in (0,1]$, from $(\ref{e1})$ we also obtain:

\begin{wn}
For every $a < x_1$, and every $u,v \in  (0,a]$, we have
$$
\Gamma (\frac{xy}{a}) \geq \Gamma (u) + \Gamma (v) - \Gamma (a).
$$
\end{wn}

\noindent
\textbf{Acknowledgement}.   The author thanks to Professor Marek W\'ojtowicz for remarks and comments which improved the quality of this paper.


\begin{thebibliography}{15}
\bibitem{Alzer1} H. Alzer, O.G. Ruehr, {\em A submultiplicative property of the psi function}, Jour. of Comput. and Applied Math. 101 (1999), 53-60.
\bibitem{Finol-Wojt} C.E. Finol, M. W\'ojtowicz, {\em Multiplicative properties of real function with applications to classical functions}, Aequationes Math. 59 (2000), 134 - 149.
\bibitem{Klamkin4} M.S. Klamkin, D.J. Newman, {\em Extensions of the Weierstrass product inequalities}, Math. Mag. 43 (1970), 137 - 140.
\bibitem{Klamkin5} M.S. Klamkin, {\em Extensions of the Weierstrass product inequalities II}, Amer. Math. Monthly 82, (1975) 741 - 742.
\bibitem{Maligranda} L. Maligranda, {\em On submultiplicativity of an N-function and its conjugate}, Aequationes math. 89, (2015), 569-573.
\bibitem{Mitrinovic} D. S. Mitrinovi\v c,  {\em Analytic Inequalities}, Springer-Verlag, Berlin 1990.
\bibitem{Nizar}  Kh. Nizar, Al-Oushoush, Moa'ath Oqielat {\em  A Generalization of Weierstrass Inequality with Some Parameters}, Australian Journal of Basic and Applied Sciences 14 (2020), 36-41.
\bibitem{Pec} J.E. Pecari\v c, M.S. Klamkin, {\em Extensions of the Weierstrass product inequalities III}, SEA Bull. Math. (II) 2 (1988), 123 - 126.
\bibitem{Pier} J. Pierpont, {\em Functions of Complex Variable}, Dover Publ., New York, 1914.
\bibitem{Stromberg} K. R. Stromberg, {\em An Introduction to Classical Real Analysis}, Brooks/Coole Publishing Company Pacific Grove, 1981.
\bibitem{WU} S. Wu,  {\em Some results on extending and sharpening the Weierstrass product inequalities }, J. Math. Anal. Appl. 308 (2005), 689-702.






\end{thebibliography}
\end{document}